\documentstyle[11pt, epsfig]{article}
\input amssym.def
\title{New Solutions of the Einstein-Dirac Equation in Dimension $n=3$.  \footnote{Supported by the SFB 288 and the Graduiertenkolleg ''Geometrie und Nichtlineare Analysis'' of the DFG.}}
\author{a short announcement by Thomas Friedrich (Berlin)}
\date{\today}

\begin{document}

\maketitle

Consider a Riemannian spin manifold of dimension $n \ge 3$ and denote by $D$ the Dirac operator acting on spinor fields. A solution of the Einstein-Dirac equation is a spinor field $\psi$ solving the equations

$$  Ric - \frac{1}{2} \, \, S \cdot g  =  \pm \frac{1}{4}  \, \, T_{\psi} \quad , \quad D(\psi) = \lambda \psi . $$

Here $S$ denotes the scalar curvature of the space, $\lambda$ is a real constant and $T_{\psi}$ is the energy-momentum tensor of the spinor field $\psi$ defined by the formula
$$ T_{\psi} (X,Y) =(X \cdot \nabla_Y \psi + Y \cdot \nabla_X \psi , \psi ) . $$

Any weak Killing spinor $\psi^*$ (WK-spinor)
$$ \nabla_X \psi^* = \frac{n}{2(n-1)} \, \, dS (X) \psi^* + \frac{2 \lambda}{(n-2)S} \, \, Ric (X) \cdot \psi^* - \frac{\lambda}{n-2} \, \, X \cdot \psi^* + \frac{1}{2(n-1)S} \, \, X \cdot dS \cdot \psi^* $$

yields a solution $\psi$ of the Einstein-Dirac equation after normalization
$$\psi = \sqrt{\frac{(n-2)\vert S \vert}{\vert \lambda \vert \vert \psi^* \vert^2}} \, \, \psi^* . $$

In fact, in dimension $n=3$ the Einstein-Dirac equation is essentially equivalent to the weak Killing equation (see [KimF]). Up to now the following 3-dimensional Riemannian manifolds admitting WK-spinors are known:
\begin{enumerate}
\item the flat torus $T^3$ with a parallel spinor;
\item the sphere $S^3$ with a Killing spinor;
\item two non-Einstein Sasakian metrics on the sphere $S^3$ admitting WK-spinors.
The scalar curvature of these two left-invariant metrics equals $S= 1 \pm \sqrt{5}$.
\end{enumerate}

The aim of this short note is to announce the existence of a one-parameter family of left-invariant metrics on $S^3$ admitting WK-spinors. This family contains the two non-Einstein Sasakian metrics with WK-spinors on $S^3$, but does not contain the standard sphere $S^3$ with Killing spinors. Moreover, any simply-connected, complete Riemannian manifold $X^3 \not= S^3$ with WK-spinors such that the eigenvalues of the Ricci tensor are constant is isometric to a space of this one-parameter family. \\
In order to formulate the result precisely, we fix real parameters $K,L,M \in {\Bbb R}$ and denote by $X^3 (K,L,M)$ the 3-dimensional, simply-connected and oriented Riemannian manifold defined by the following structure equations:
$$ \omega_{12} = K \sigma^3 \quad , \quad \omega_{13}= L \sigma^2 \quad , \quad \omega_{23} = M \sigma^1  , $$

or, equivalently:

$$ d \sigma^1 =(L-K) \sigma^2 \wedge \sigma^3 \quad , \quad d \sigma^2 = (M+K) \sigma^1 \wedge \sigma^3 \quad , \quad d \sigma^3 = (L-M) \sigma^1 \wedge \sigma^2 . $$

The 1-forms $\sigma^1, \sigma^2 , \sigma^3$ are the dual forms of an orthonormal frame of vector fields. Using this frame the Ricci tensor of $X^3 (K,L,M)$ is given by the matrix
$$ Ric = \left( \begin{array}{ccc}
-2KL & 0 & 0\\
0 & 2KM & 0 \\
0 & 0 & -2LM
\end{array} \right) . $$

{\bf Theorem:} {\it Let $X^3 \not= S^3$ be a complete, simply-connected Riemannian manifold such that:
\begin{itemize}
\item[a)] the eigenvalues of the Ricci tensor are constant;
\item[b)] the scalar curvature $S \not= 0$ does not vanish.
\end{itemize}

If $X^3$ admits a WK-spinor, then $X^3$ is isometric to $X^3 (K, L,M)$ and the parameters are a solution of the equation 
$$ - K^2 L(L-M)^2 M+ L^3 M^3 + KL^2 M^2 (M-L) + K^3 (L-M)(L+M)^2=0 \hspace{1cm} (*)$$

Conversely, any space $X^3 (K,L,M)$ such that  $(K,L,M) \not= (0,0,0)$ is a solution of  $(*)$ admits two WK-spinors for one and only one WK-number $\lambda$. With respect to the fixed orientation of $X^3 (K,L,M)$ we have the two cases:
$$ \lambda = + \frac{S}{2 \sqrt{2}} \, \, \sqrt{ \frac{S}{S^2 - |Ric|^2}} \quad \quad \mbox{if $-K <M$}$$
$$ \lambda = - \frac{S}{2 \sqrt{2}} \, \, \sqrt{ \frac{S}{S^2 - |Ric|^2}} \quad \quad  \mbox{if $M < -K$} . $$

\bigskip

The spaces $X^3 (K,L,M)$ are isometric to $S^3$ equipped with a left-invariant metric.\\}

\bigskip

{\bf Remark:} If the parameters $K=M$ coincide, the solution of the equation $(*)$ is given by 
$$ L = \frac{1}{4} K(1 - \sqrt{5}) \quad , \quad L= \frac{1}{4} K(1+ \sqrt{5})  $$

and we obtain the Ricci tensors

$$ Ric = \left( \begin{array}{ccc}
\frac{1}{2} K^2 (\sqrt{5}-1) & 0 & 0\\
0 & 2K^2 & 0\\
0 & 0 & \frac{1}{2} K^2 (\sqrt{5}-1) \end{array} \right)$$

or

$$ Ric = \left( \begin{array}{ccc}
- \frac{1}{2} K^2 (1+ \sqrt{5}) & 0 & 0\\
0 & 2K^2 & 0\\
0 & 0 & - \frac{1}{2} K^2 (1+ \sqrt{5}) \end{array} \right) . $$

\bigskip

The non-Einstein-Sasakian metrics on $S^3$ occur for the parameter $K=1$ (see [KimF]).\\

\newfont{\graf}{eufm10}
\newcommand{\sod}{\mbox{\graf so}(3)}

{\bf Remark:} Using the standard basis of the Lie algebra $\sod$ we can write the left-invariant metric of the space $X^3 (K,L,M)$ in the following way:
$$ \left( \begin{array}{ccc}
\frac{1}{|M-L||K+M|} & 0 & 0\\
\\
0 & \frac{1}{|K-L||M-L|} & 0\\
\\
0 & 0 & \frac{1}{|K-L||K+M|} \end{array} \right) . $$

The equation $(*)$ is a homogeneous equation of order six. The transformation $(K,L,M) \to (\mu K, \mu L, \mu M)$ corresponds to a homothety of the metric. Therefore - up to a homothety - the moduli space of solutions is a subset of the real projective space ${\Bbb P}^2 ({\Bbb R})$ given by the equation $(*)$. This subset is a configuration of six curves in ${\Bbb P}^2 ({\Bbb R})$ connecting  the three points $[K:L:M] = [1:0:0], [0:1:0], [0:0:1]$, corresponding to flat metrics. 

\begin{center}

\[
\epsfig{figure=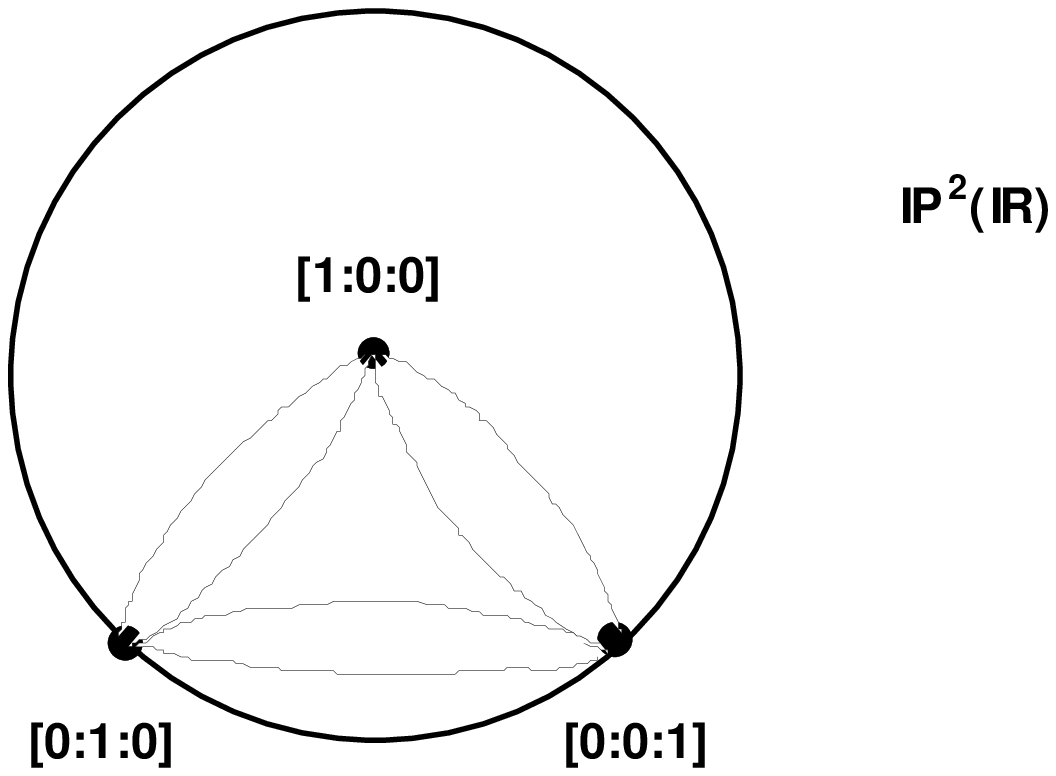,width=8cm}
\]

\end{center}

In particular, we have constructed two paths of solutions of the Einstein-Dirac equation deforming the non-Einstein Sasakian metrics on $S^3$.\\

\bigskip

The proof of the Theorem as well as the complete computations will be published in a furthercoming paper of the author (see [F]).\\

\end{document}